# LAPLACE TRANSFORM EFFECTIVENESS IN THE $M|G|\infty$ QUEUE BUSY PERIOD PROBABILISTIC STUDY

FERREIRA Manuel Alberto M. (PT)

**Abstract.** The Laplace transform is a widely used tool in the study of probability distributions, often allowing for a probability density functions and distribution functions simpler determination and being a "moments generating function". In this paper it is considered a situation not so simple, as it is the case of the M|G|∞ queue busy period length distribution. Attention will also be given the respective tail Laplace transform. Then, in the context of an open queues network, which nodes behave as M|G|∞ queues, the Laplace transform will be used to construct an algorithm to determine the Laplace transform of the global service time length of a customer during their stay on the network distribution.

## 1  Introduction

In the $M|G|\infty$ queue, customers arrive according to a Poisson process at rate $\lambda$, upon its arrival receive immediately a service with time length *d. f.* $G(\cdot)$ and mean $\alpha$. The traffic intensity is $\rho=\lambda\alpha$, see for instance [13].

As it happens for any queue, in the $M|G|\infty$ queue activity there is a sequence of idle and busy periods. For this queue the study of the busy period length distribution is very important since, as it is part of its definition, a customer must find immediately an available server upon its arrival. So, it is important for the manager to know how many, and how long, servers must be in prevention, see [4, 10 and 20].

In the next two sections, the busy period length and the busy period length tail Laplace transforms will be studied, by this order, in some of their most important details.

The $M|G|\infty$ queue busy period length distribution is called *B*, the *d. f. B(t)* and the *p. d. f. b(t)*. The Laplace transform is denoted $\bar{B}(s)$.

A **network of queues** is a collection of nodes, arbitrarily connected by arcs, across which the customers travel instantaneously and

 - There is an arrival process associated to each node,

- There is a **commutation process** which commands the various costumers' paths.

The arrival processes may be composed of **exogenous arrivals**, from the outside of the collection, and of **endogenous arrivals**, from the other collection nodes. Call

$$\Lambda = \begin{bmatrix} \lambda_1 \\ \lambda_2 \\ \vdots \\ \lambda_J \end{bmatrix} \quad (1.1)$$

the network exogenous arrival rates vector, where $\lambda_j$ is the exogenous arrival rate at node $j$ and $\lambda = \sum_{j=1}^{J} \lambda_j$.

A network is **open** if any customer may enter or leave it. A network is **closed** if it has a fixed number of customers that travel from node to node and there are neither arrivals from the outside of the collection nor departures. A network open for some customers and closed for others is said **mixed**.

The commutation process rules, for each costumer that abandons a node, which node it can visit then or if it leaves the network. In a network with $J$ nodes, the matrix

$$P = \begin{bmatrix} p_{11} & p_{12} & \cdots & p_{1j} \\ p_{21} & p_{22} & \cdots & p_{2j} \\ \vdots & \vdots & \vdots & \vdots \\ p_{j1} & p_{j2} & \cdots & p_{jj} \end{bmatrix} \quad (1.2)$$

is the commutation process matrix, being $p_{jl}$ the probability of a customer, after ending its service at node $j$, go to node $l$, $j, l = 1, 2, \ldots J$. The probability $q_j = 1 - \sum_{l=1}^{J} p_{jl}$ is the probability that a customer leaves the network from node $j, j = 1, 2, \ldots J$. Call now $\gamma_j$, the total – from the outside of the network and from the other nodes – customers arrival rate at node $j$ and

$$\Gamma = \begin{bmatrix} \gamma_1 \\ \gamma_2 \\ \vdots \\ \gamma_J \end{bmatrix} \quad (1.3)$$

the network exogenous arrival rates vector. If the network is stable, the following equality – **traffic equations** – holds:

$$\Gamma^T = \Lambda^T + \Gamma^T P \quad (1.4).$$

Note that they may be written as $\Gamma^T = \Lambda^T (I - P)^{-1}$. For more details on networks of queues see [2 and 24].

In section 4, for open networks of queues, which nodes are $M|G|\infty$ queues, it will be constructed an algorithm to determine the Laplace transform of the distribution of the global service time length of a client during their stay on the network, see [1].

This work finishes with the presentation of a conclusions section and a short list of references.

## 2  The $M|G|\infty$ Queue Busy Period Length Laplace Transform

The $M|G|\infty$, queue busy period length Laplace transform is, see [13, 15 and 24],

$$\bar{B}(s) = 1 + \lambda^{-1}\left(s - \frac{1}{\int_0^\infty e^{-st - \lambda \int_0^t [1-G(v)]dv} dt}\right) \quad (2.1).$$

From (2.1) it is easy to obtain the following expression for the mean busy period length, see [1, 23 and 24]:

$$E[B] = \frac{e^\rho - 1}{\lambda} \quad (2.2),$$

for any service time distribution.

Inverting (2.1) it is achieved the $M|G|\infty$ Queue Busy Period Length $p.\,d.\,f.$, see [3]:

$$b(t) = G(0)\delta(t) + (1 - G(0))\left[\frac{d}{dt}\left(e^{-\lambda \int_0^t [1-G(v)]dv} \frac{G(t) - G(0)}{1 - G(0)}\right)\right]$$

$$* \sum_{n=0}^{\infty} \left[\frac{d}{dt}(1 - e^{-\lambda \int_0^t [1-G(v)]dv})\right]^{*n} \quad (2.3),$$

where $\delta$ is the Dirac delta and $*$ is designated convolution operator. For constant service time with value $\alpha$, expression (2.3) becomes, see [3],

$$b(t) = \sum_{n=0}^{\infty} g(t) * \left[\frac{dA(t)}{dt}\right]^{*n} e^{-\rho}(1 - e^{-\rho})^n \quad (2.4)$$

where $g(t) = \frac{dG(t)}{dt}$ and $A(t) = \begin{cases} \frac{1-e^{-\lambda t}}{1-e^{-\rho}}, t < \alpha \\ 1, t \geq \alpha \end{cases}$.

For a $M|G|\infty$, queue, if the service time $d.\,f.$ belongs to the collection, see [13, 14],

$$G(t) = 1 - \frac{1}{\lambda} \frac{(1-e^{-\rho})e^{-\lambda t - \int_0^t \beta(u)du}}{\int_0^\infty e^{-\lambda w - \int_0^w \beta(u)du} dw - (1-e^{-\rho})\int_0^t e^{-\lambda w - \int_0^w \beta(u)du} dw}, t \geq 0, -\lambda \leq \frac{\int_0^t \beta(u)du}{t}$$

$$\leq \frac{\lambda}{e^\rho - 1} \qquad (2.5)$$

the busy period length *d. f.* (achieved inverting $\frac{1}{s}\bar{B}(s)$) is

$$B(t) = \left(1 - (1-G(0))\right)\left(e^{-\lambda t - \int_0^t \beta(u)du} + \lambda \int_0^t e^{-\lambda w - \int_0^w \beta(u)du} dw\right)$$

$$* \sum_{n=0}^\infty \lambda^n (1-G(0))^n \left(e^{-\lambda t - \int_0^t \beta(u)du}\right)^{*n}, t \geq 0, -\lambda \leq \frac{\int_0^t \beta(u)du}{t}$$

$$\leq \frac{\lambda}{e^\rho - 1} \qquad (2.6).$$

**Notes:**

-The demonstration may be seen in [13],

-For $\frac{\int_0^t \beta(u)du}{t} = -\lambda$, $G(t) = B(t) = 1, t \geq 0$ in (2.5) and (2.6), respectively,

-For $\frac{\int_0^t \beta(u)du}{t} = \frac{\lambda}{e^\rho - 1}$, $B(t) = 1 - e^{-\frac{\lambda}{e^\rho - 1}t}, t \geq 0$, only exponential in (2.6),

-If $\beta(t) = \beta$ (constant)

$$G(t) = 1 - \frac{(1-e^{-\rho})(\lambda+\beta)}{\lambda e^{-\rho}(e^{(\lambda+\beta)t}-1)+\lambda}, t \geq 0, -\lambda \leq \beta \leq \frac{\lambda}{e^\rho - 1} \qquad (2.7)$$

and

$$B^\beta(t) = 1 - \frac{\lambda + \beta}{\lambda}(1-e^{-\rho})e^{-e^{-\rho}(\lambda+\beta)t}, t \geq 0, -\lambda \leq \beta \leq \frac{\lambda}{e^\rho - 1} \qquad (2.8),$$

a mixture of a degenerate distribution at the origin and an exponential distribution.∎

The expression (2.1) is equivalent to $(\bar{B}(s) - 1)C(s) = \lambda^{-1}sC(s)$ where $C(s) = \int_0^\infty e^{-st - \lambda \int_0^t [1-G(v)]dv} \lambda(1 - G(t))dt$. Differentiating n times to *s* and using Leibnitz´s formula, it is achieved the expression, see [4, 13 and 15],

$$E[B^n] = (-1)^{n+1}\left\{\frac{e^\rho}{\lambda}nC^{(n-1)}(0) - e^\rho \sum_{p=1}^{n-1}(-1)^{n-p}\binom{n}{p}E[B^{n-p}]C^{(p)}(0)\right\}, n = 1,2,\ldots \quad (2.9)$$

$$C^{(n)}(0) = \int_0^\infty (-t)^n e^{-\lambda \int_0^t [1-G(v)]dv} \lambda(1-G(t))dt, n = 0,1,\ldots$$

that permits to write exact formulae for the $E[B^n], n = 1,2,\ldots$ through a recurrent process. The efficiency of these formulae depends on the possibility of calculating the various $C^{(n)}(0)$. But, for instance, for the distribution given by (2.8) it is obtained

$$E[B^n] = \frac{\lambda+\beta}{\lambda}(1-e^{-\rho})\frac{n!}{(e^{-\rho}(\lambda+\beta))^n}, n = 1,2,\ldots \quad (2.10).$$

And in the situation of constant service times

$$C^{(0)}(0) = 1-e^{-\rho} \text{ and } C^{(n)}(0) = -e^{-\rho}(-\alpha)^n - \frac{nC^{(n-1)}(0)}{\lambda}, n = 1,2,\ldots \quad (2.11),$$

so, the calculations are quite simple.

## 3 The $M|G|\infty$ queue busy period length tail Laplace transform

Call $H(t) = 1 - B(t), t \geq 0$ the $M|G|\infty$ busy period length tail and $\bar{H}(s)$ it's Laplace Transform. The essential result in this section is, see [16]:

**Lemma 3.1**

$$1 - G(t) = \lambda^{-1} \frac{L^{-1}\left[\frac{1}{\lambda\bar{H}(\cdot)+1}\right](t)}{\int_0^t L^{-1}\left[\frac{1}{\lambda\bar{H}(\cdot)+1}\right](v)dv}, t \geq 0 \quad (3.1).\blacksquare$$

**Notes**:

-The symbol $L^{-1}$ means inverse Laplace transform,

-Obviously $1 - G(t)$ is the service length tail,

-Expression (3.1) allows determining the service tail corresponding to a busy period tail. But there are situations for which $1 - G(t)$ does not match to a tail. So, it is fundamental to look for conditions that $\bar{H}(s)$ must satisfy to guarantee that a tail is obtained through expression (3.1).∎

Using Polya's Theorem, see for instance [22], and defining

$$a(t) = \frac{\lambda\bar{H}(-it)}{e^\rho - 1}, i = \sqrt{-1} \quad (3.2)$$

the following result was achieved, see again [16]:

**Lemma 3.2**

If $a(t)$ is a real function different from $\frac{1}{1-e^{\rho}}$ and satisfy the conditions:

- $\dfrac{\frac{d^2 a(t)}{dt^2}\left(a(t)+\frac{1}{e^{\rho}-1}\right)-2\left(\frac{da(t)}{dt}\right)^2}{a(t)-\frac{1}{e^{\rho}-1}} > 0, t > 0$

- $\lim\limits_{t\to\infty} a(t) = 0$

$1 - G(t)$ obtained through expression (3.1) is a tail. ∎

**Note**:

-Distributions that do not fulfil this lemma cannot be $M|G|\infty$, busy period length distributions. ∎

## 4      An algorithm to compute the global service time distribution in an open network of $M|G|\infty$ queues through Laplace transforms

An open network of queues with infinite servers in each node, with Poisson process exogenous arrivals, may be looked like a $M|G|\infty$ queue. The service time is the sojourn time of a customer in the network, see [17 and 19].
Note that the sojourn time is the mixture of the sums of the services corresponding to each path that a customer may have in the network. The total time spent in a path by a customer distribution is so the convolution of the service time distributions in each node belonging to the path, since those service times are independent. Each one of these convolutions is a parcel in the mixture which weight is given by the path probability. Each path starts in a node $j$ with probability $\frac{\lambda_j}{\lambda}$ and ends in node $k$ with probability $1 - \sum_{j=1}^{J} p_{kj}$.

As the Laplace transform of a convolution of two functions is the product the two those functions Laplace transforms and having in mind the traffic equations seen above (expression (1.4)):

- Denote $S$ the sojourn time of a costumer in the network and $S_j$ its service time at node $j, j = 1,2,...J$. Be $G(t)$ and $G_j(t)$ the $S$ and $S_j$ distribution functions, respectively, and $\bar{G}(s)$ and $\bar{G}_j(s)$ the Laplace transforms,
- Define

$$\Lambda(s) = \begin{bmatrix} \lambda_1 \bar{G}_1(s) \\ \lambda_2 \bar{G}_2(s) \\ \vdots \\ \lambda_J \bar{G}_J(s) \end{bmatrix} \text{ and } P(s) = \begin{bmatrix} p_{11}\bar{G}_1(s) & p_{12}\bar{G}_2(s) & \cdots & p_{1J}\bar{G}_J(s) \\ p_{21}\bar{G}_1(s) & p_{22}\bar{G}_2(s) & \cdots & p_{2J}\bar{G}_J(s) \\ \vdots & \vdots & & \vdots \\ p_{J1}\bar{G}_1(s) & p_{J2}\bar{G}_2(s) & & p_{JJ}\bar{G}_J(s) \end{bmatrix} \quad (4.1)$$

- It results

$$\bar{G}(s) = \sum_{n=0}^{\infty} (\lambda^{-1}\Lambda^T(s)P^n(s)(I-P)A) \quad (4.2),$$

- And finally, using the Leontief's matrix properties, the global service time distribution Laplace transform service time for a customer during its permanence in the network is given by

$$\bar{G}(s) = \lambda^{-1}\Lambda^T(s)\bigl(I-P(s)\bigr)^{-1}(I-P)A \quad (4.3),$$

where $I$ is the identity matrix with the same order as $P$ and $A$ is a column with J 1`s, for the Laplace Transform service time, confer with [17].

So, the problem in terms of Laplace transforms is operationally simple. The "problems" arrive when inverting the Laplace transform. As usual the situation is not bad when using exponential expressions.

## 5    Conclusions

In this text the operational qualities and also the research incentive of the Laplace transform is well known. The results presented are of both types: purely quantitative and of theoretical scope.
In studies about stochastic processes, of which queues are part, it is very common to use this tool. The situation is like that which occurs in differential equations: great simplification in operational matters, not always accompanied by comparable simplification when it is necessary to reverse the transform.
It happens that in the case of stochastic processes, it is often possible to collect the fruits of the research without recourse to inversion.
From the presented results it is worth noting the formula (4.3) for its simplicity and evident utility, where the qualities of the Laplace transform are quite explored.

**Acknowledgement**

This work is financed by national funds through FCT - Fundação para a Ciência e Tecnologia, I.P., under the project UID/MULTI/4466/2016. Furthermore, I would like to thank the Instituto Universitário de Lisboa and ISTAR-IUL for their support.

**Current address**

**Manuel Alberto M. Ferreira, Professor Catedrático**
INSTITUTO UNIVERSITÁRIO DE LISBOA (ISCTE-IUL)
BRU – IUL, ISTAR-IUL
Av. das Forças Armadas, 1649-026 Lisboa, Portugal
Tel.: + 351 21 790 37 03. FAX: + 351 21 790 39 41,
E-mail: manuel.ferreira@iscte.pt